\documentclass[10pt]{article}
\usepackage{epsfig}
\usepackage{amsmath}
\usepackage{amssymb}
\usepackage{amscd}
\usepackage{tabularx}
\usepackage{float}
\usepackage{pstricks}
\usepackage{pst-grad}
\usepackage{float}
\usepackage{theorem}
\usepackage{enumerate}
\newpsobject{showgrid}{psgrid}{subgriddiv=1,griddots=10,gridlabels=6pt}
\def\eqnum#1{}

\begin{document}
\title {Small values of the maximum for the integral of fractional Brownian motion}
\author
{
G. Molchan \thanks{Observatoire de la Cote d'Azur, CNRS UMR 6529, BPU 229, 06304, Nice Cedex 4,
 France, molchan@mitp.ru}\footnotemark[1] \footnotemark[2] 
\and 
A. Khokhlov  \thanks{International Institute of Earthquake Prediction Theory
and Mathematical Geophysics 79, b2,  Warshavskoe shosse 113556 Moscow, Russia.
khokhlov@mitp.ru, khokhlov@ipgp.jussieu.fr} 
}
\date{}
\maketitle
\begin{abstract}
We consider the integral of fractional Brownian motion (IFBM)
and its functionals $\xi_T$ on the intervals $(0,T)$ and $(-T,T)$
of the following types:
the maximum $M_T$, the position of the maximum,
the occupation time above zero etc. We show how the
asymptotics of $P(\xi_T<1)=p_T, T\to \infty$, is related
to the Hausdorff dimension of Lagrangian regular points
for the inviscid Burgers equation with FBM initial velocity.
We produce computational evidence in favor of a
power asymptotics for $p_T$. The data do not reject the
hypothesis that the exponent $\theta$
of the power law is related to the similarity parameter $H$ of
fractional Brownian motion as follows: $\theta =-(1-H)$  for the
interval (-T,T) and $\theta =-H(1-H)$ for $(0,T)$.
The point $0$ is special in that IFBM and its
derivative both vanish there.

\end{abstract}

\newpage

\bigskip

{\bf 1. Introduction}

\bigskip
Sinai [18] and Frisch and associates [17] initiated in 1992 the study of
fractal and multifractal properties of solutions of the inviscid Burgers
equation with initial velocity $u_0(x)$ specified by a self-similar random
process. That last circumstance guarantees that the solution is
self-similar in the large. In particular, one could be interested in
finding the Hausdorff dimension of the set of regular Lagrangian points S
that describe the initial locations of those fluid particles which have not
collided until a fixed time $t_0$. The original model of $u_0(x)$
was fractional Brownian motion (FBM), $b_H(x)$, with similarity
parameter $0<H<1$.

By now the Sinai-Frisch program has been carried out for special
Markovian models of $u_0(x)$ alone: Sinai [18] has found the
dimension $S$ for
Brownian motion case, i.e., $u_0(x)=b_{1/2}(x)$; Bertoin [3] discovered
for this case that the solution
$u(t=t_0,x)$ admits of an exact probabilistic
description in terms of a stable Levy process. One can then find a
multifractal description of the solution $x\to u(t_0,x)$ (the relevant
references are [6, 7, 8]).
Additionally, Bertoin [3] found that the Hausdorff dimension of Lagrangian
regular points is $h$, if $u_0(x)$ is a stable L\'evy process of index
$\alpha =h^{-1}\in (1,2]$ with no positive jumps (see also [20]).

The nonmarkovian case $u_0(x)=b_H(x)$, $H\ne 1/2$ has proved extremely
difficult for analysis. Handa [5] found simple arguments to derive
a lower bound on the dimension of $S$, namely, $\mbox {\rm dim}S\ge H$.
The exact equality $\mbox {\rm dim}S=H$ is known as a
hypothesis [17, 19] since 1992. Among methods developed for analyzing the
nonmarkovian case $u_0(x)$, the Sinai approach is of particular interest.
For the case $u_0(x)=b_{1/2}(x)$ this method [18] relates the estimation
of the dimension of $S$ to the asymptotic behavior of the probability

\begin{eqnarray}
p_T = P\{ \xi (x)<1, \quad x \in \Delta_T \}
\eqnum{1}
\end{eqnarray}
for integral Brownian motion

\begin{eqnarray}
\xi (x) = \int^x_0b_{1/2}(s)ds \quad \mbox {\rm and} \quad
\Delta_T = (0,T),\quad T \gg 1.
\eqnum{2}
\end{eqnarray}

As a matter of fact (see below), one has to deal with a problem
that is rather popular in physical and technical applications: find the
probability of a large excursion for a random process
$\eta (x)$, i.e., $P\{ \eta (x)>0, 1<x<T\}$, $T \gg 1$.
A review of the problem can be found in [14].
Sinai has shown that the quantity $p_T\cdot T^{1/4}$ is bounded away
from 0 and $\infty$ as $T\to \infty$ under the conditions (2).
That estimate was repeatedly
refined  and generalized [9, 10, 11].

We show below that the upper bound $\mbox {\rm dim}S\le H$ under
the conditions $u_0(x)=b_H(x)$ follows from an estimate of
$p_T$ for the integral of fractional Brownian motion (IFBM):
$\xi (x) = \int^x_0b_H(s)ds$ when
considered in the bilaterally expanding interval $\Delta_T=(-T,T)$.

The work [15] clarifies the asymptotic problem
of $p_T$ for  intervals $(0,T)$ and $(-T,T)$  in the case
of fractional Brownian motion:
$\xi (x)=b_H(x)$. It transpires that in this case

\begin{eqnarray}
\ln p_T = -(1-H)\ln T (1+o(1)), \quad \Delta_T = (0,T).
\nonumber
\end{eqnarray}
On the other hand, when
$\Delta_T=\{ x: |x|<T\}$, the leading term in the log
asymptotics of $p_T$ is independent of $H$. More generally, suppose
$b_H(x), x\in R^d$ is FBM with multidimensional time; in that case

\begin{eqnarray}
\ln p_T = -d\ln T (1+o(1)), \quad \Delta_T = \{ x\in R^d:|x|<T\}.
\nonumber
\end{eqnarray}

The last asymptotics is due to the fact that the probability density for
the position of the maximum of FBM exists in the sphere $\{ |x|<1\}$.
A generalization of this fact is given below.

We present theoretical and computational evidence in favor of the
following asymptotics for IFBM:

\begin{eqnarray}
\ln p_T = \left \{ \begin{array}{ll}
-(1 - H)\ln T (1 + o(1)), \quad \quad \Delta_T = (-T,T) \\
-H(1 - H)\ln T (1 + o(1)), \quad \Delta_T = (0,T).
\end{array}\right.
\nonumber
\end{eqnarray}

The first of these asymptotic expressions corroborates the
hypothesis $\mbox {\rm dim}S=H$, so is not unexpected, while the second is,
considering that the exponent $\theta (H)=H(1-H)$
has the point of symmetry $H=1/2$.

Because IFBM is a self-similar process, the distribution of its
maximum in $\Delta =(0,1)$ or (-1,1), $F_{\max}(x)$, is related to
$p_T$ through $p_T=F_{\max}(T^{-(1+H)})$.
Importantly, our calculation was performed for a series of statistics:
the maximum \,\,   $M=\max\limits_{\Delta}\mbox {\rm IFBM}$;\,\,
the position of the maximum  $M$,  $|G|$;  the  occupation
time    $A^+=\int_\Delta {\bf 1}_{\xi(x)>0}dx$    of IFBM  above zero; and
the rightmost zero of IFBM in $(0,T)$, $Z$. The distributions of these
statistics (one should use $F_{\max}(x^{1+H})$
when $M$ is considered) have identical asymptotics as $x\to 0$, but
depend on interval type: $\Delta =(0,1)$ or (-1,1). When
$\Delta =(-1,1)$, they provide independent evidence in favor of the
hypothesis $\mbox {\rm dim}S=H$.

The rest of this paper is organized as follows. Section 2 reduces
the evaluation of $\mbox {\rm dim}S\le H$ to the asymptotic
distributions of $M$, $G$, $A^+$ and $Z$
near zero. Section 3 discusses the modeling
of IFBM, while section 4 presents numerical evaluations of the
distributions listed above and some theoretical arguments to support our
conclusions.

\bigskip
{\bf 2. Regular Lagrangian points and the nonexceedance of level}

\bigskip
We now define more exactly the notions used in Introduction.
We consider the Burgers equation

\begin{eqnarray}
\partial_tu + u\partial_xu = \nu u_{xx}, \quad \nu \downarrow 0,
\eqnum{3}
\end{eqnarray}
with continuous initial conditions $u(0,x)=u_0(x)$ and the velocity
potential $U(x)=\int^x_0u_0(x)dx=o(x^2)$, $x\to \infty$.
The solution at $t_0=1$ has the form $u(x)=x-a(x)$, where $a(x)$ can
be found from $U(x)$ as follows. Construct a convex minorant $C(x)$
for $U(x)+x^2/2$. In
that case its derivative $C'(x)$ is nondecreasing and has finite limits
from the left and from the right. We now complete the definition of $C'(x)$
in continuity on the right. In that case, according
to Hopf (see, e.g., [18, 21]),
$a(x)$ is identical with the inverse function of $C'(x)$. The set of
points where $C'(x)$ is increasing, i.e., the topological support of the
measure  $dC'(x)$ or the closure of the set $\{ a(x), x\in R\}$,
defines the set $S$ of regular Langrangian
points in the Burgers problem. The dynamics of completely inelastic
particles on $R^1$ can be related to the Burgers equation:
each  infinitesimal particle located at $x$
has a mass $dx$ and an initial momentum $dU(x)$. On colliding the
particles coalesce and continue movement following the
conservation laws of mass and momentum. Particles that have not collided
until time $t_0=1$ make up the set $S$ in the Lagrangian coordinates. The
initial conditions $u_0(x)$ will be considered to be fractional Brownian
motion $b_H(x)$, i.e., a Gaussian process with zero mean and structural
function $E|b_H(x)-b_H(y)|^2=|x-y|^{2H}$ where $0<H<1$. In virtue of
the Kolmogorov theorem the paths of $b_H(x)$ can be
treated as continuous a.s. The process $b_H(x)$ is self-similar, i.e.,
$b_H(\Lambda x)\stackrel {\rm d}{=}\Lambda^Hb_H(x)$,
where $\stackrel {\rm d}{=}$ denotes equality of finite-dimensional
distributions.

{\bf Theorem 1.} 1. The set of regular Lagrangian points $S$ in the
Burgers problem (3) with $u_0(x)=b_H(x)$ has a.s. dimension $H$, if
for any $\varepsilon >0$ and
$T\to \infty$ one of the following requirements is fulfilled:

\begin{eqnarray}
\nonumber
&(A)& \quad P\{ y(x): = \int^x_0b_H(s)ds < 1, \,\,
x \in \Delta_T\} <  T^{-(1-H)+\varepsilon}, \\
\nonumber
&(B)& \quad P(y(x) < 0, \,\,
x \in \Delta_T, \, |x| > 1) <  T^{-(1-H)+\varepsilon}, \\
\nonumber
&(C)& \quad P(|G(\Delta_T)| < 1) <  T^{-(1-H)+\varepsilon}, \\
\nonumber
&(D)& \quad P\{ \int_{\Delta_T} {\bf 1}_{y(x)>0}dx < 1,\,\,
|G(\Delta_T)| < T\} <  T^{-(1-H)+\varepsilon},
\nonumber
\end{eqnarray}
where $\Delta_T=(-T,T)$, $G(\Delta_T)$ is the position of the maximum
of $y(x)$ in $\Delta_T$.

2. If one of type $A-D$ probabilities $p_T$ has an asymptotics of the
form $\log p_T=-\theta \log T(1+o(1))$, the probabilities of the other
types have the same asymptotics. This statement also holds
for $\Delta_T=(0,T)$ with the probability
$P(Z_T<1)$  in addition to $(A-D)$, where $Z_T$ is the rightmost zero
of $y(x)$ in $(0,T)$.

The proof of the theorem will be preceded by two lemmas.

{\bf Lemma 1.} $\mbox {\rm dim}S\le H$, if
for any $\varepsilon >0$  there exists a $\delta_0=\delta_0(\varepsilon)$
such that one has for arbitrary $x\in R^1$:

\begin{eqnarray}
P(S\cap B(x,\delta)\ne \phi) < \delta^{(1-H)-\varepsilon}, \,\,
\delta < \delta_0
\eqnum{4}
\end{eqnarray}
where $B(x,\delta)$ is a ball of radius $\delta$ centered at $x$.

{\it Proof.} Cover the interval $\Delta =[a,b]$ with intervals
$B_i(\delta)$ of length $\delta$ with overlappings of length
$\delta /2$. Consider the measure $\mu (dx)=dC'(x)$ with support $S$.
The elements $\tilde {B}_i$ in $\{ B_i(\delta)\}$ for
which $\mu (B_i)>0$ will then form a cover $S\cap \Delta$. In view of (4)

\begin{eqnarray}
\nonumber
E\sum |\tilde {B}_i(\delta)|^{H+2\varepsilon}=
E\sum |B_i(\delta)|^{H+2\varepsilon}{\bf 1}_{\mu(B_i)>0}  \\
< \delta^{H+2\varepsilon}\cdot 2|\Delta|\delta^{-1}\cdot
\delta^{(1-H)-\varepsilon} = c\delta^\varepsilon \qquad
\nonumber
\end{eqnarray}
where $|\Delta|$ is the length of $\Delta$. By Chebyshev's inequality

\begin{eqnarray}
P(\sum |\tilde {B}_i|^{H+2\varepsilon} > a) < c\delta^{\varepsilon} /a.
\nonumber
\end{eqnarray}
Consider a sequence $\delta_n$ such that
$\sum \delta_n^\varepsilon < \infty$. The Borel-Cantelli
lemma then yields

\begin{eqnarray}
\sum |\tilde {B}_i(\delta_n)|^{H+2\varepsilon} < a,
\quad n > n(\omega)
\nonumber
\end{eqnarray}

Since $a$ is arbitrary:

\begin{eqnarray}
\lim\sup_n\sum |\tilde {B}_i(\delta_n)|^{H+2\varepsilon} = 0
\quad \mbox {\rm a.s.}
\nonumber
\end{eqnarray}

However, in that case one has
$\mbox {\rm dim}(S\cap \Delta) \le H+2\varepsilon$. Since
$\varepsilon > 0$ and $\Delta$ are arbitrary, one has
$\mbox {\rm dim}S\le H$.

{\bf Lemma 2.} The conditions of Lemma 1 are fulfilled, if

\begin{eqnarray}
P\{ \int^x_0b_H(s)ds < 1, \quad |x| < T\} < T^{-(1-H)+\varepsilon}, \quad
\forall \varepsilon > 0.
\nonumber
\end{eqnarray}
as $T \to \infty$.

{\it Proof.} The process
$y(x)=\int^x_0b_H(s)ds+x^2/2$  can be represented in the form

\begin{eqnarray}
\nonumber
y(x)&=&\int^c_0(b_H(s)+s)ds + (b_H(c)+c)(x-c) + \\
    &+&\int^x_c[b_H(s)-b_H(c)+(s-c)]ds = L(x') +
\int^{x'}_0(\tilde{b}_H(s)+s)ds,
\nonumber
\end{eqnarray}
where $L(x')$ is a linear function of $x'=x-c$, and
$\tilde{b}_H(x)=b_H(c+x)-b_H(c) \stackrel{\rm d}{=}b_H(x)$.
The convex minorants of $y$ and
$\tilde{y}=\int^{x'}_0(\tilde{b}_H(s)+s)ds$
differ by the linear function $L(x')$. Hence the fractal properties of the
measure $\mu(dx)=dC'(x)$ are invariant under
translation along the $x$-axis (this observation is due to U. Frisch).
Consequently, it is sufficient to  prove (4)  for
$S'=S\cap (-\delta /2,\delta /2)$.

Let $ \Delta =(-\delta /2,\delta /2)$ contain a point of growth
$x_0$ for the measure $d\mu$. That means that the curve
$f(x)=U(x)+x^2/2$ and its convex
minorant $C(x)$ do not lie below the tangent of $f(x)$ at the point $x_0$,
and $C(x_0)=f(x_0)$. The event $\{ x_0\in \Delta \}$, to be called
$A$ here, can be written as

\begin{eqnarray}
\nonumber
A&=&\{ \exists x_0: |x_0|<\delta /2;\,\,\int^x_0(b_H(s)+s)ds  \\
 &>&\int^{x_0}_0(b_H(s)+s)ds + (b_H(x_0)+x_0)(x-x_0),\,\,
\forall x\in R^1\}.
\nonumber
\end{eqnarray}

Let us modify event $A$ to become $A_1$, i.e., we assume that the equality
in the formulation of $A$ is true for $|x|<1$ only. To emphasize
the fact that $A_1$ depends on the process $b_H(x)+x=\xi (x)$, we will
write $A_1=A_1[\xi]$.

One has

\begin{eqnarray}
P(A)\le P(A_1) = E\,{\bf 1}_{A_1[b_H+\varphi]} =
E\,{\bf 1}_{A_1[\tilde{b}_H]}\pi (\tilde{b}_H),
\eqnum{5}
\end{eqnarray}
where $\varphi (x)=x$ and $\pi$ is the Radon-Nikodim derivative of two
Gaussian measures corresponding to the processes $\tilde{b}_H-\varphi$
and $\tilde{b}_H$ in [-1,1]. Note that $\tilde{b}_H$ is an FBM process.
The function $\varphi$ is smooth and vanishes at zero.
For this reason the above measures are
mutually absolutely continuous [16]. By the Cameron-Martin relation
$\ln \pi (\tilde{b}_H)$ is a Gaussian variable with mean
$-c^2_H/2$ and variance $c_H^2$,
where $c_H=\Vert \varphi \Vert$ and $\Vert \cdot \Vert$ is
the  norm in Hilbert space $H_B$
of functions on $\Delta =[-1,1]$
with reproducing kernel $B(x,y)=Eb_H(x)b_H(y)$. The
constant $c_H$ is finite and can be found in explicit form as
indicated by Molchan and Golosov [16].

Applying H\"older's inequality to the right-hand side of (5), one gets

\begin{eqnarray}
P(A) <  P(A_1[\tilde{b}_H])^{1-\varepsilon}
(E\pi^{1/\varepsilon})^\varepsilon =
P(A_1[\tilde{b}_H])^{1-\varepsilon}c_\varepsilon,
\eqnum{6}
\end{eqnarray}
where $c_\varepsilon =\exp(\frac {1}{2} (\varepsilon^{-1}-1)c^2_H)$.

We now evaluate $P(A_1[b_H])$. One has

\begin{eqnarray}
\nonumber
P(A_1[b_H])&=&P\{ \exists x_0:|x_0|<\delta /2;   \\
\nonumber
\int^x_0b_H(s)ds&>&\int^{x_0}_0b_H(s)ds + b_H(x_0)(x-x_0), \,|x|<1\}\\
\nonumber
&=&P\{ \exists x_0:|x_0|<1/2; \int^x_0b_H(s)ds > a(x_0) + b(x_0)x, \, |x|<T\} ,
\nonumber
\end{eqnarray}
where

\begin{eqnarray}
\nonumber
T=\delta^{-1}, \quad
|a(x_0)|&=&|\int^{x_0}_0b_H(s)ds - b_H(x_0)x_0| < 2M, \\
\nonumber
|b(x_0)|&=&|b_H(x_0)| < M = \max_{|x|<1/2}|b_H(x)|.
\nonumber
\end{eqnarray}

We will use the Fernique inequality [4]

\begin{eqnarray}
P(M > \bar{c}_Hu) < \exp(-u^2/2) = T^{-a}, \quad u>u_0,
\nonumber
\end{eqnarray}
where $u=u_T=\sqrt{2a\ln T}$, $c_H$ being a constant;
the value of $a$ will be chosen later on.
From this it follows that

\begin{eqnarray}
\nonumber
P(A_1[b_H])&<&P\{ A_1[b_H], \,\ M < \bar{c}_Hu_T\} + T^{-a} \\
\nonumber
&<&P\{ \int^x_0b_H(s)ds > -2\bar{c}_Hu_T - \bar{c}_Hu_T|x|,
\,\, |x|<T\} + T^{-a} \\
\nonumber
&=&P\{ \int^x_0b_H(s)ds < u_T\bar{c}_H(2+|x|), \, |x|\le T\} + T^{-a} \\
&=&P\{ \int^x_0b_H(s)ds < 4\lambda_T^{-1} + 2|x|, \, |x|<T'\} + T^{-a},
\eqnum{7}
\end{eqnarray}
where
$T'=T/\lambda_T$, \, $u_T\bar{c}_H=2\lambda_T^H$,
$\lambda_T=\mbox {\rm const} \cdot (\ln T)^{1/2H}$.
Here we have used the fact that $b_H(x)$ is a
self-similar process and modified the interval $|x|\le T$ to
become $|x| < T'$.

Define the function

\begin{eqnarray}
\varphi_1(x) = 2x\,{\bf 1}_{|x|<1} + 2\mbox {\rm sgn}(x){\bf 1}_{|x|>1} =
\frac{2}{\pi i} \int[e^{ix\lambda} - 1]\frac{\sin \lambda}{\lambda^2}
\,d\lambda.
\eqnum{8}
\end{eqnarray}
In that case (7) can be continued to get

\begin{eqnarray}
P(A_1[b_H]) \le
P\{ \int^x_0(b_H(s) - \varphi_1(s))ds < F(x),  \,\, |x|<T'\} + T^{-a},
\eqnum{9}
\end{eqnarray}
where

\begin{eqnarray}
F(x) = \left \{ \begin{array}{ll}
-x^2 + 2|x| + 4\lambda_T^{-1}, \quad |x|<1, \\
1 + 4\lambda_T^{-1}, \qquad \qquad \quad \, |x|>1.
\end{array}\right.
\nonumber
\end{eqnarray}
When $T$ is large, one has $F(x)<2$. For this reason the last estimate
will merely become less precise, when $F$ is replaced with $F(x)=2$.
The right-hand side of (9) can be evaluated by repeating the steps that
have led to (5, 6). The substitution of $\tilde{b}_H$ for
$b_H-\varphi_1$ combined with H\"older's inequality yield

\begin{eqnarray}
P(A_1[b_H]) <
P\{ \int^x_0b_H(s)ds < 2,  \,\,
|x|<T'\}^{1-\varepsilon}c^1_\varepsilon + T^{-a},
\nonumber
\end{eqnarray}
where $c^1_\varepsilon = \exp((\varepsilon^{-1}-1)b^2_H/2)$,
$b^2_H = \Vert \varphi_1\Vert^2_T \le \Vert \varphi_1\Vert^2_\infty$.
Here  $\Vert \cdot \Vert_T$ is the norm on $H_B$ for the interval
$(-T,T)$.
The spectral representations of the kernel

\begin{eqnarray}
B(t,s) = Eb_H(t)b_H(s) = k^{-1}_H
\int(e^{ix\lambda} - 1) (e^{-ix\lambda} - 1) |\lambda|^{-1-2H}\,d\lambda
\nonumber
\end{eqnarray}
and $\varphi_1$ (see (8)) yield

\begin{eqnarray}
\Vert \varphi_1\Vert^2_\infty = k_H\int \bigg\vert \frac{2\sin \lambda}
{\pi \lambda^2}\bigg\vert^2|\lambda|^{1+2H}\,d\lambda < \infty ,
\nonumber
\end{eqnarray}
where $k_H=\int |e^{i\lambda}-1|^2|\lambda|^{-1-2H} d\lambda$.

The final result is

\begin{eqnarray}
P(A) \le P(A_1[b_H])^{1-\varepsilon}c_\varepsilon <
(p_{T'}^{1-\varepsilon}\, c^1_\varepsilon + T^{-a})
^{1-\varepsilon}c_\varepsilon,
\eqnum{10}
\end{eqnarray}
where $p_{T'} = P\{ \int^x_0b_H(s)ds < 2,\,|x|<T'\}$.

Let $p_T<T^{-(1-H)+\varepsilon_1}$ for large $T$. Take $a>1-H$ and
choose $\varepsilon $ from the requirement
$c^1_\varepsilon \cdot c_\varepsilon =T'^{\varepsilon_1}$,
i.e., $\varepsilon =c\varepsilon^{-1}_1/\ln T'$.
Inequality (10) can then be continued:

\begin{eqnarray}
P(A) < c_1T'^{-(1-H)+2\varepsilon_1},
\nonumber
\end{eqnarray}
where $c_1=\exp (2(1-H)c\varepsilon_1^{-1})$. Recalling that
$T'=T(\ln T)^{-1/2H}\cdot c_2$, one obtains the desired estimate
$P(A)<T^{-(1-H)+3\varepsilon_1}$, $T>T_0(\varepsilon_1)$, $T=\delta^{-1}$.

{\it Proof of Theorem 1}. The inequality $\mbox {\rm dim}S\ge H$ was
derived by Handa [5]. The opposite inequality
$\mbox {\rm dim}S\le H$  follows from Lemmas 1 and 2 and
condition $(A)$ of Theorem 1. To prove the theorem under condition $(B)$,
we note that the event
$\{ \int^x_0b_H(s)ds < c,\, |x|<T\}$
can be represented as

\begin{eqnarray}
\{ \int^x_0(b_H(s) - \varphi (s))ds < \psi (x), \quad |x|<T\},
\nonumber
\end{eqnarray}
where $\varphi, \psi$ are smooth finite functions: $\psi \equiv 0$
when $|x|\ge 1$ and $\psi >0$ when $|x|<1$, while $\varphi =0$ when
$|x|<1/2$ and $|x|>1$. Repeating the translation procedure for the
samples: $b_H(s)-\varphi (s)\to \tilde{b}_H(s)$ and using H\"older's
inequality, we get

\begin{eqnarray}
\nonumber
p_T:&=&P\{ \int^x_0b_H(s)ds < 1,\,\, |x|<T\} \\
\nonumber
&<& c_\varepsilon P\{ \int^x_0b_H(s)ds < \psi (x),\,\,|x|<T\}^{1-\varepsilon} \\
\nonumber
&<&c_\varepsilon P\{ \int^x_0b_H(s)ds < 0,\,\,1<|x|<T\}^{1-\varepsilon},
\nonumber
\end{eqnarray}
where $c_\varepsilon =\exp (\frac{1}{2}\varepsilon^{-1}\cdot c^2_\varphi)$,
$c_\varphi <k_H\int|\hat{\varphi}(\lambda)|^2|\lambda|^{1+2H}$  and
$\hat{\varphi}$ is the Fourier transform
of $\varphi$. One has $c_\varphi <\infty$, because $\varphi$ is smooth
and finite.
Choose $\varepsilon =\varepsilon_T$ from the requirement
$c_\varepsilon ={{{\cal L}}}_T$, where ${{{\cal L}}}_T$ is a slowly
varying function. Take ${{{\cal L}}}_T=\ln T$, say, then
$\varepsilon_T^{-1}=c\ln \ln T$. The result is

\begin{eqnarray}
p_T < {{{\cal L}}}_TP(\int^x_0b_H(s)ds < 0,\,\,1<|x|<T)^{1-\varepsilon_T}.
\nonumber
\end{eqnarray}

When $(B)$ holds, one has

\begin{eqnarray}
p_T < (T^{-(1-H)+\varepsilon_1})^{1-\varepsilon_T}{{{\cal L}}}_T <
T^{-(1-H)+\varepsilon_2}, \,\, T\gg 1,
\nonumber
\end{eqnarray}
i.e., the implication $(B)\to (A)$ is true.

The inequality

\begin{eqnarray}
P(\int^x_0b_H(s)ds < 0,\,\,1<|x|<T) < P(|G(-T,T)|<1),
\nonumber
\end{eqnarray}
where $G(\Delta)$ is the position of the maximum of IFBM in $\Delta$,
yields the implication $(C)\to (B)$.
Lastly,  under $(B)$ the position of the maximum of IFBM is $|G|<1$,
while the occupation time of IFBM above 0 is below 2. Hence $(D)\to (B)$.

Let us prove $(A)\to (C)$. Below, $G_T$ is the position
of the maximum of IFBM in $(-T,T)$ and
$M_a=\max\limits_{|x|<a}\mbox {\rm IFBM}$. One has

\begin{eqnarray}
P(|G_T| < 1) < P(|G_T| < 1), \,\, M_1 < c_T) + P(M_1>c_T).
\nonumber
\end{eqnarray}

If $c_T=\sqrt{2a\ln T}$, then the Fernique estimate [4] yields

\begin{eqnarray}
P(M_1 > c_T) < cT^{-a'},\,\,T > T_0,
\nonumber
\end{eqnarray}
where $a'=a/\sigma^2$ and
$\sigma^2=E|\mbox {\rm IFBM}(1)|^2=(2H+2)^{-1}$. Also,

\begin{eqnarray}
P(|G_T| < 1, \,M_1 < c_T) < P(M_T < c_T).
\nonumber
\end{eqnarray}

Since IFBM is self-similar, one has
$M_T\stackrel{\rm d}{=}\lambda^{1+H}M_{T'}$, when $T=\lambda T'$.
Take $\lambda$ from the requirement $\lambda^{1+H}=c_T$. Then
$P(M_T < c_T) = P(M_{T'} < 1)$.

\noindent To sum up,

\begin{eqnarray}
P(|G_T| < 1) < P(M_{T'} < 1) + o(T^{-a}),
\nonumber
\end{eqnarray}
where $T'=cT(\ln T)^{-\rho}$, $\rho =(1-H)^{-1}$, while the parameter
$a>0$ is arbitrary. When $a>(1-H)$,
the implication (A)$\to$(C) is obvious.

We are going to prove (A)$\to$(D). Let $A^+_T$ be the occupation
time of $y(x)=\mbox {\rm IFBM}(x)$ above zero in $\Delta_T=(-T,T)$. One has

\begin{eqnarray}
P(A^+_T < 1,\,|G_T| < T) \le P(M_T < c_T) + P(M_T > c_T,\,
A^+_T < 1,\,|G_T| < T),
\nonumber
\end{eqnarray}
where $c_T$ will be specified below.

Let $\Delta_T=\cup \Delta_k$,\,\,$\Delta_k=(k,k+1)$ and
$M_k=\max \{ y(x), x\in \Delta_k\}$. If the event
${\cal B}=\{ M_T > c_T,\,A^+_T < 1,\,|G_T| < T\}$ occurs, one will
have the following for the interval $\Delta_k$ which contains
$G_T$:  $M_k>c_T,\, y(x)$ and
$y'(x)=b_H(x)$ have zeroes in $\Delta_k$. Indeed,
if $y(x)\ne 0$, then $y(x)>0$ in $\Delta_k$ and $A^+_T\ge 1$.
Consequently,

\begin{eqnarray}
P({\cal B}) < \sum_k P\{ \max ((y(x_1) - y(x_2)) > c_T,\,\,
x_1,x_2\in \Delta_k),\,\, S_k\}: = \sum_kp_k,
\nonumber
\end{eqnarray}
where $S_k$ means that $b_H(x)$ has a zero in $\Delta_k$.

We are going to evaluate $p_k$:

\begin{eqnarray}
\nonumber
p_k&<&P\{ \max[(y(x_1) - y(x_2)), \, x_1,x_2\in \Delta_k] > c_T,\,\,
|b_H(k)| < c_T/2\} \\
   &+& P\{ |b_H(k)| > c_T/2,\,\, S_k\} : = p_{k,1} + p_{k,2}.
\nonumber
\end{eqnarray}

One has

\begin{eqnarray}
\nonumber
p_{k,2}&<&P\{ \max(|b_H(x_1) - b_H(x_2)|, \, x_1,x_2\in \Delta_k) >
c_T/2\} \\
       &=&P\{ \max(|b_H(x_1) - b_H(x_2)|, \, x_1,x_2\in \Delta_0) >
c_T/2\}.
\nonumber
\end{eqnarray}
Here we have used the fact that $b_H(x)$ has stationary increments.
In virtue of the Fernique inequality [4]

\begin{eqnarray}
p_{k,2} < c\, \exp (-\frac{1}{2} (c_T/c_b)^2),
\nonumber
\end{eqnarray}
where $c$ is an absolute constant, while $c_b$ is a function of $H$.

One proceeds in a similar manner to evaluate $p_{k,1}$:

\begin{eqnarray}
y(x_1) - y(x_2) = \int^{x_1}_{x_2}(b_H(s) - b_H(k))dx + b_H(k)(x_1 - x_2).
\nonumber
\end{eqnarray}
If $\max [y(x_1)-y(x_2)]>c_T$  in $\Delta_k\times \Delta_k$ and
$|b_H(k)|<c_T/2$, then

\begin{eqnarray}
\max \int^{x_1}_{x_2}[b_H(s) - b_H(k)]ds > c_T/2.
\nonumber
\end{eqnarray}
Consequently,

\begin{eqnarray}
\nonumber
p_{k,1}&<&P\{ \max[\int^{x_1}_{x_2}(b_H(s) - b_H(k))ds, \,
x_1,x_2\in \Delta_k] > c_T/2\} \\
       &=&P\{ \max [\int^{x_2}_{x_1}b_H(s)ds, \, x_1,x_2\in \Delta_0] >
c_T/2\}.
\nonumber
\end{eqnarray}
Here again, we have used the relation
$b_H(x)-b_H(k)\stackrel{\rm d}{=}b_H(x-k)$ with a fixed $k$.
The use of the Fernique inequality [4] yields

\begin{eqnarray}
p_{k,1} < c\, \exp (-\frac{1}{2} (c_T/2c_y)^2)
\nonumber
\end{eqnarray}
where $c_y$ is a function of $H$. Combining the estimates of $p_{k,1}$
and $p_{k,2}$ and assuming $c_T=\max (c_b,2c_y)\sqrt{2a \log T}$, one gets

\begin{eqnarray}
p_k = p_{k,1} + p_{k,2} < cT^{-a}.
\nonumber
\end{eqnarray}
However, one then has $P({\cal B}) < 2cT^{-a+1}$
and

\begin{eqnarray}
\nonumber
P(A^+_T < 1, \, |G_T| < T)&\le &P(M_T < c_T) + O(T^{-a+1}) \\
&=&P(M_{T'} < 1) + O(T^{-a+1}),
\nonumber
\end{eqnarray}
where $T'=c T(\log T)^{-\rho}$, $\rho =(1-H)^{-1}$. Hence $(A)\to (D)$.

Consider the second part of the theorem.
Let $p_T({\Theta})$ be the probabilities that appear in the first part of
the theorem, where ${\Theta}$ denotes the condition $A, B, C$ or $D$.
It has been shown above that, when $T\gg 1$,

\begin{eqnarray}
\nonumber
p_T(A)&<&p_T(B)^{1-\varepsilon_T}{{{\cal L}}}_T, \,\,p_T(B) < p_T(C), \\
p_T(C)&<&p_{T'}(A) + O(T^{-a}), \, \,p_T(D) < p_{T'}(A) + O(T^{-a+1}),
\nonumber
\end{eqnarray}
where $a>0$ is any fixed number, ${\cal L}_T$ is a slowly varying function
of $T$, $\varepsilon_T=o(1)$, $T\to \infty$ and
$T'=c T(\ln T)^{-\rho}$, $\rho =(1-H)^{-1}$.
A trivial corollary of these is that all
the $p_T({\Theta})$ have the asymptotics
$\ln p_T({\Theta})=-\theta \ln T(1+o(1))$, provided the asymptotics
is true for at least a single quantity of the type ${\Theta}=A, C$ or $D$.

Our proof has not relied significantly on the type of the
interval $\Delta_T$: $(-T,T)$ or $(0,T)$. For this reason our conclusion
that the asymptotics of $p_T({\Theta})$ are identical also holds for $(0,T)$.

We conclude by noting that, if $\Delta_T=(0,T)$, then
$p_T(B) = P(y(x) > 0,\, 1<t<T)$.
Consequently, if $Z_T$ is the rightmost zero of $y(x)$ in $(0,T)$, then
$P(Z_T<1) = 2p_T(B)$.

\bigskip

{\bf 3. The generation of IFBM}

\bigskip

We are going to use Monte Carlo techniques in order to evaluate
the probabilities $p_T({\Theta})$ with ${\Theta}=A,C,D$ in Theorem 1
for the process $y(x)=\int^x_0b_H(s)ds$
in the following intervals of $\Delta_T: (0,T)$ and $(-T,T)$.
The probabilities in question are small, $p_T\to 0$ as $T\uparrow \infty$,
hence the IFBM generation should be exact for a discrete sequence
$\{ x_k, k=1,...,T\}$.
Since $y(x)$ is a self-similar process, it is sufficient to use integer
points $x_k=k$. In that case
$\{ y(k/T)\} \stackrel{\rm d}{=}\{ T^{-(1+H)}y(k)\}$, while the
probabilities $p_T({\Theta})$ can obviously be expressed in terms of the
statistics $M=\max\limits_{\Delta_1}y(x)$,\,
$G=\arg\max\limits_{\Delta_1}y(x)$,\,\,$A^+=\int_{\Delta_1}{\bf 1}_{y>0}dx$\,
of the process $\{ y(x), x\in \Delta_1\}$, where $\Delta_1=(0,1)$
or $(-1,1)$, as follows:

\begin{eqnarray}
p_T(A)=F_M(T^{-(1-H)});\quad p_T(C)=F_{|G|}(T^{-1}); \quad
p_T(D)=\hat{F}_A(T^{-1})F_G(1-0)
\nonumber
\end{eqnarray}
where $F_\xi$ is the distribution of $\xi$ and $\hat{F}_A$ is
the conditional distribution of $A^+$ given $|G|\ne 1$.

{\it The generation of} $\{ y(k), k=0,...,T\}$. The sequence
$\{ y(k), k=0,...,T\}$ is Gaussian and has
stationary second increments, i.e., the sequence

\begin{eqnarray}
\eta_k = y(k-1) - 2y(k) + y(k+1),\quad k=1,...,T-1,
\eqnum{11}
\end{eqnarray}
has a Toeplitz correlation matrix $[\mu_{i-j}]$, where

\begin{eqnarray}
\mu_k = c_{q}[|k-2|^{q} - 4|k-1|^{q} + 6|k|^{q} - 4|k+1|^{q} +
|k+2|^{q}]
\eqnum{12}
\end{eqnarray}
and $c_{q}=[2q (q -1)]^{-1}$, $q =2H+2$.

The second differences (11) combined with the initial
conditions $y(0)=0$ and $y(1)$ are sufficient to uniquely reconstruct
the sequence $\{y(k), k=0,...,T\}$. One can assign $y(1)$ by using
the decomposion  $y(1) = \hat{y}(1) + y^\perp(1)$
into the predictable part
$\hat{y}(1)=E(y(1)|\eta_1,...,\eta_{T-1})$ of $y(1)$
and the part $y^{\perp}(1)$ that
cannot be predicted from the data $\{ \eta_i,  i=1,...,T-1)\}$. In that case

\begin{eqnarray}
y(1) = \sum^{T-1}_{k=1}z_k\eta_k + \sigma \varepsilon_0
\eqnum{13}
\end{eqnarray}
Here, ${\bf z}=(z_1,...,z_{T-1})'$ is the solution of the linear equation:

\begin{eqnarray}
[\mu_{i-j}]^{T-1}_1{\bf z} = {\bf m},
\eqnum{14}
\end{eqnarray}
where the vector ${\bf m}$ has the components

\begin{eqnarray}
m_k = Ey(1)\eta_k = {\bf \Delta }[q |k|^{q -1} - |k|^q + |k-1|^q]c_q
\nonumber
\end{eqnarray}
and ${\bf \Delta}$ is the difference operator of second order:
${\bf \Delta }f(k)=f(k-1)-2f(k)+f(k+1)$. The second
term is $\sigma \varepsilon_0=y^\perp (1)$, where $\varepsilon_0$ is
the standard Gaussian random
variable which is independent of $\{ \eta_1,...,\eta_{T-1}\}$;

\begin{eqnarray}
\sigma^2 =E[y^\perp (1)]^2 = q^{-1} - \sum_{1\le k<T}z_km_k,
\nonumber
\end{eqnarray}
because $q^{-1}=E|y(1)|^2$ and
$E|\hat{y}(1)|^2=\sum^{T-1}_1z_km_k$.

It thus appears that the exact generation of the sequence
$\{y(k), k=0,...,T\}$
reduces to the generation of the stationary Gaussian sequence
$\{ \eta_k,  k=1,...,T-1)\}$ with
correlation function (12) and to the solution of the linear equation (14).

{\it The generation of} $\{ y(k),\,|k|<T/2\}$. For generating $y(x)$ in
a bilateral interval, we note the following. Assume that $y(x)$ is
IFBM in $(0,T)$, while $\tilde{y}(x)$
is IFBM in $(-T',T''), T'+T''=T$. In that case

\begin{eqnarray}
\{ y(x) - y(T') - y'(T')(x-T'),\,x\in(0,T)\}
\stackrel{\rm d}{=}\{ \tilde{y}(x-T'), \,x\in (0,T)\}
\nonumber
\end{eqnarray}

The left-hand side provides a key to how one is to transform the
sequence $\{ y(k),\,k=0,...,T\}$ into an IFBM sequence that starts
from the point $0<k_0<T$. To do this one must also find the derivative
$y'(k_0)$. In a similar way as above:

\begin{eqnarray}
y'(k_0) = E\{y'(k_0)|\eta_1,...,\eta_{T-1}\} + E\{y'(k_0)|
\varepsilon_0\} + y'^{\perp}(k_0)
\nonumber
\end{eqnarray}
where the first two terms correspond to the predictable
part of $y'(k_0)$ based on the data
$\{ \eta_1,...,\eta_{T-1}, \varepsilon_0\}$, while the third term
corresponds to  the unpredictable part of $y'(k_0)$. The predictable
part is

\begin{eqnarray}
E\{y'(k_0)|\eta_1,...,\eta_{T-1}, \varepsilon_0\} =
\sum^{T-1}_{k=1}z'_k\eta_k + a\varepsilon_0,
\nonumber
\end{eqnarray}
where $(z'_1,...,z'_{T-1})$ is the solution of (14) with the right-hand
side ${\bf m}'=(m'_1,...,m'_{T-1})$.
The components of ${\bf m}'$ are

\begin{eqnarray}
m'_k = E y'(k_0)\eta_k = {\bf \Delta }[|k|^{q -1} + |k_0-k|^{q -1}
\mbox {\rm sgn} (k_0-k)]qc_q.
\nonumber
\end{eqnarray}

One has $a=E y'(k_0)\varepsilon_0$. From (13) one derives
$E y(1) y'(k_0)=\sigma a+\sum^{T-1}_{k=1}z_km'_k$.
Hence

\begin{eqnarray}
\sigma a = c_q \cdot q [(q -1)k_0^{q- 2}+1-k_0^{q -1}+
(k_0-1)^{q -1}] - \sum^{T-1}_{k=1}z_km'_k.
\nonumber
\end{eqnarray}

One has
$y'^{\perp}(k_0) = \sigma' \varepsilon'$
where $\varepsilon'$ is a standard Gaussian variable that is independent
of $(\eta_1,...,\eta_{n-1},\varepsilon)$. The variance of the
unpredictable part $\sigma'^2$ can be found from the relation

\begin{eqnarray}
k_0^{2H} = E [y'(k_0)]^2 =
\sum^{T-1}_{k=1}z'_km'_k + a^2 + \sigma'^2.
\nonumber
\end{eqnarray}

To sum up, the exact generation of
$\{ y(k),  k=0,...,T; y(k_0)=y'(k_0)=0\}$
requires that an equation like (14) should be solved twice.

{\it The generation of} $\{ \eta_k\}$. Bardet et al. [2] provide
a review of the methods which allow generation of Gaussian
stationary sequences with a prescribed correlation function. We
use the progressive Schur algorithm [1], which is a
Levinson-Durbin method. The Generalized
Schur algorithm can be used in the framework of this method for
fast solution of equations like (14) by the Gohberg-Semenkul
formula [1]. The generation of $\{ \eta_k,  k=1,...,T-1\}$
by this method requires $O(T^2)$ operations. The computation
is organized so as to minimize the amount of calculation needed for
generating $N$ IFBM samples; the computational complexity is a linear
function of $N$ and the storage capacity is of order $O(T^2)$.

The parameters $T$ and $N$ are equally important in the problem
discussed. However, it appears that the increase of $T$ would
not be effective from {\it a priori} considerations. The argument is
as follows. We are interested in the distributions of $M$, $G$, $A^+$ near
zero where they are expected to behave like $x^\theta{\cal L}(x)$,
where ${\cal L}$ is a slowly
varying function. Since the approximation to IFBM is discrete, these
distributions contain a positive atom at zero of size
$p_0(T) = P(y(k) \le 0,\, \, k\in \Delta_T)$.
The probability is doubled for the statistic $Z$ (the rightmost zero
of $y(t), t\in (0,T)$).

Note that
$p_0(T) \ge P(y(x),\,\, |x|>1, \, \, x\in \Delta_T).$ Therefore,
if the bounds given by Theorem 1 are explicit, $p_0(T)$ should
not decrease faster than  $T^{-\theta} {{{\cal L}}}'(T)$.
The rate of decrease becomes very low, when $\theta <1$.
The expected value is $\theta =1-H$ for the interval $(-T,T)$
and $\theta \le 1/4$ for $(0,T)$ (see below).

Given the above situation, it only remains to increase $N$ alone.
This will provide good accuracy for empirical estimates of the
distributions of $M$, $G$, $A^+$, $Z$ for discrete time.
The discrete-time distributions cannot converge to the basic continuous
ones faster that $T^{-\theta (H)}$ that we have at $x=0$.

\bigskip

{\bf 4. Evaluation of $\theta (H)$ and related results}

\bigskip

{\it The bilateral IFBM process}.

Figure 1 presents estimates of the distribution of $|G_{1/2}|$
($G_a$ is the position of the maximum of
$y(x)=\int^x_0b_H(s)ds$ in the interval $(-a,a)$).
The distributions are given for $H: 0.1-0.9$ at
increments of 0.1. The estimates are based on $N=50,000$ samples
of $y(x)$ with sample size $T=8194$. We use $|G_{1/2}|$
to demonstrate the asymptotics of $p_T$ from Theorem 1, because
in this case the distribution discontinuity at $x=0$
bends the graph of the distribution on a log-log scale near zero
to a lesser extent. We recall that the discontinuity
of size $p_0(T)$ occurs for all the statistics $M$, $G$ and $A^+$,
and is due to the presence of negative excursions in $y(x)$ on a
discrete grid of $x$. The curves in Fig. 1
are well consistent with the asymptotics of type
$F_G(x)=x^\theta {{{\cal L}}}(x)$, where ${{{\cal L}}}$ is a slowly
varying function and $\theta =1-H$. Assuming
${\cal L}(x)=\mbox {\rm constant}$ in $(x_-,x_+)$, one can construct
the estimate of maximum likelihood $\hat{\theta}$ for $\theta$.

\begin{figure}[h]
\centerline{
\begin{pspicture}(0,0)(10,5)
\psset{dimen=inner}
\rput(5,2.5){
\epsfig{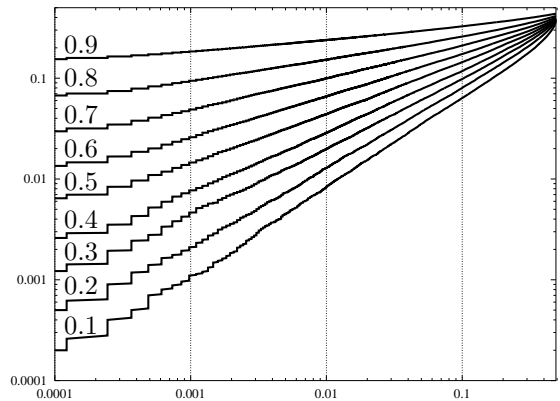}
}
\rput(2.3, 0.87){$0.1$}
\rput(2.3, 1.4){$0.2$}
\rput(2.3, 1.85){$0.3$}
\rput(2.3, 2.275){$0.4$}
\rput(2.3, 2.78){$0.5$}
\rput(2.3, 3.22){$0.6$}
\rput(2.3, 3.68){$0.7$}
\rput(2.3, 4.16){$0.8$}
\rput(2.3, 4.61){$0.9$}
%%\showgrid
\end{pspicture}
}
\caption  {Distributions of $|G|$ for the position
of maximum $G$ of the IFBM process in interval ($-\frac1{2},\frac1{2}$). Various values of $H$
are shown at the left hand part of each curve.
}
\end{figure}

The estimate $\hat{\theta}$ was computed in intervals of the
form $(10^{-3}, 10^{-2})\times i$, $i=1-5$. The choice of the
initial point 0.001 is related to the fact that the sampling
interval over time is $\Delta =T^{-1}=0.00012$,
so that all curves in  $(\Delta, 10\Delta)$  slightly change
their slopes due to the discontinuity of discrete-time distributions
at zero. The deviation of  $\hat{\theta}$  from the hypothetical
$\theta_0=1-H$ does not exceed 0.03, i.e., is less than 6\% when
$H\le 0.5$. The error is large ($\ge 10\%$) for $0.6\le H\le 0.9$
because of smallness of $\theta_0$. For $H\ge 0.6$ we have
a slow convergence of the discrete-time
distributions with $\Delta =T^{-1}\downarrow 0$  and difficulties in the
choice of small $x$ to estimate $\theta$.

The argument to be given below provides a partial
explanation of the nature of the asymptotics

\begin{eqnarray}
P(|G_T|<1) = T^{-(1-H)}{{{\cal L}}}_T.
\eqnum{15}
\end{eqnarray}

The function $y(x)$ is differentiable, hence the position of the global
maximum, $G_T$, belongs to the zero set of $b_H(x)$ or to
the end-points of $(-T,T)$. For this reason it should seem that the
local time
$l(x)=\lim\limits_{\varepsilon \to 0} \frac{1}{2\varepsilon}
\int^x_0{\bf 1}_{|b_H(s)|<\varepsilon}ds$
is the natural time scale in our problem of the
maximum of $y(x)$, i.e., it is more natural to study
$\tilde{y}(l)=y(x(l))$ instead of $y(x)$, where $x(l)$ is
the inverse function of $l(x)$ which is continuous
on the right. The process was first treated by Vergassola et al. [19]
and independently
used by Isozaki and Watanabe [11] to prove the Sinai
asymptotics for $H=1/2$. It is a known fact [12] that $l(x)$ is
a continuous self-similar process with parameter
$h=1-H$. Consequently, $l(T)=O(T^{1-H})$, and (15) means, roughly
speaking, that $P(|\tilde{G}_L|<1)=L^{-1}{{{\cal L}}}_L$ where
$\tilde{G}_L$ is the location of the maximum of
$\tilde{y}(l)$ in $(-L,L)$. A more exact statement can be made.
Let ${{{\cal L}}}_{T,i}$ be slowly varying functions that decrease
as $T\to \infty$.

{\bf Statement 2.} If\quad  a) $G_T$ and $l(\pm T)$ are weakly
dependent, i.e.

\begin{eqnarray}
P\{ |G_T| < 1, \, \,  |l(\pm T)| > T^{1-H}{\cal L}_{T,1})\} >
P(|G_T| < 1){{{\cal L}}}_{T,2}
\eqnum{16}
\end{eqnarray}

\noindent and
$\mbox {\rm b)} \quad
P(|\tilde{G}_L| < 1) < [L{{{\cal L}}}_L]^{-1}$,\,
then $P(|G_T|<1)<[T^{(1-H)}{\cal L}_{T,3}]^{-1}$.

{\it Proof}. By (16) we have

\begin{eqnarray}
\nonumber
P(|G_T| < 1)&<&{\cal L}_{T,2}^{-1}\,P\{ |G_T| < 1,\,
|l(\pm T)| > T^{1-H}{{{\cal L}}}_{T,1}, \\
|l(\pm G_T)|&<& (a \ln T)^H\} + {\cal L}_{T,2}^{-1}\,
P\{ |l(\pm 1)| > (a \ln T)^H\}
\eqnum{17}
\end{eqnarray}
The second term on the right-hand side admits of an upper bound
according to [22]: ${{{\cal L}}}_{T,2}^{-1}\cdot T^{-ca}$
where $c$ is an absolute constant. It follows that (17) can be continued:

\begin{eqnarray}
\le {\cal L}_{T,2}^{-1}\,\,[P\{ |\tilde{G}_{L(T)}| < (a \ln T)^H\}
+ O(T^{-ac})],
\nonumber
\end{eqnarray}
where $L(T)=T^{1-H}{{{\cal L}}}_{T,1}$.

The process $y(l)$ is self-similar; therefore,

\begin{eqnarray}
P\{ |\tilde{G}_{L(T)}| < (a \ln T)^H)\} = P\{ |\tilde{G}_{L^*}| < 1\}
< [L^*{\cal L}_{L^*}]^{-1},
\nonumber
\end{eqnarray}
where $L^*=L(T)(a \ln T)^{-H}$. Combining the resulting estimates
and bearing in mind that the constant $"a"$ is arbitrary, one gets
Statement 2. $\diamond$

Note that the condition (b) of Statement 2
is automatically satisfied for $H=1/2$, because
$\tilde{y}$ is a stable Levy process. The next theorem  shows that
this condition also holds for general self-similar processes with
stationary increments (SSSi). It is the case for $H=1/2$ again.

{\bf Theorem 3.} Let $\xi(t)$, $\xi(0)=0$ be an SSSi-process
for which sample paths have only discontinuities of the first
kind, a.s. Let $M$ be $\sup \{ \xi (s), s\in [0,1] \}$ and $G$
be the leftmost position of $M$:

\begin{eqnarray}
G = \inf \{ t\in [0,1]: \sup_{s\to t}\xi (s) = M\}
\nonumber
\end{eqnarray}

Then

1) $G$ has
a continuous probability density $\psi (t)$ in (0,1), and

\begin{eqnarray}
\psi (t) \le \psi (s) \max \biggr (\frac {s}{t}, \frac {1-s}{1-t}\biggl ),
\quad \forall t,s\in (0,1),
\nonumber
\end{eqnarray}
i.e., $\psi \equiv 0$ or $\psi >0$.

2) if $\psi \not\equiv 0$, then the position $G_T$ of the supremum of
$\xi (t)$ in $(-Ta, T(1-a))$
satisfies the following relation:

\begin{eqnarray}
P(|G_T|<1) = \frac {\psi (a)}{2T} (1+o(1)),\quad T\to \infty.
\eqnum{18}
\end{eqnarray}

{\it Remark}. Theorem 3 shows that the asymptotics of type (18) is due
to the presence of a nonzero distribution density for the position
of the supremum of the SSSi process in (0,1). That fact was first
pointed out in [15] for the process $b_H(x)$. The proof of the general
case is nearly identical with that given in [15], so it is relegated
to the Appendix.

{\it The unilateral IFBM process.}

\begin{figure}[h]
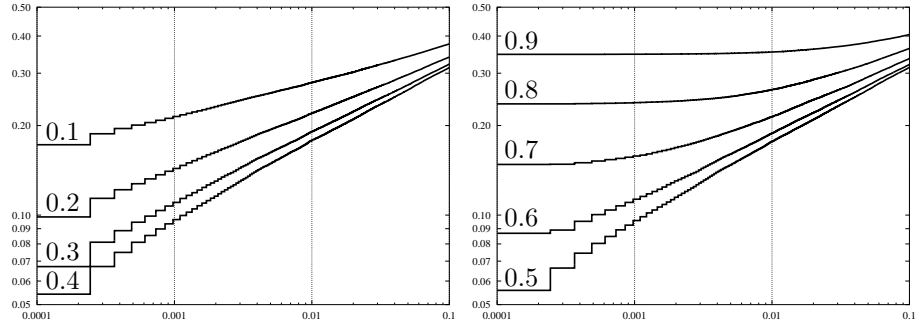

\begin{pspicture}(0,0)(12,4)
\psset{dimen=inner}
\rput(6,2.){
\begin{minipage}{0.49\linewidth}
\centerline{ 
\epsfig{file=./Fig_IFBM_minus.epsi, angle=270, width=0.99\linewidth}
}
\end{minipage}
\begin{minipage}{0.49\linewidth}
\centerline{ \epsfig{file=./Fig_IFBM_plus.epsi, angle=270, width=0.99\linewidth}}
\end{minipage}
}
\rput(0.7, 2.4){$0.1$}
\rput(0.7, 1.45){$0.2$}
\rput(0.7, 0.8){$0.3$}
\rput(0.7, 0.41){$0.4$}
\rput(6.8, 0.45){$0.5$}
\rput(6.8, 1.25){$0.6$}
\rput(6.8, 2.15){$0.7$}
\rput(6.8, 2.95){$0.8$}
\rput(6.8, 3.6){$0.9$}
%%\showgrid
\end{pspicture}
\caption  {Distribution of the position
of maximum $G$ for the IFBM process. 
Shown are the parts of  curves related to the  interval ($10^{-4},10^{-1}$)  together with the corresponding values of $H$}
\label{BNG2}
\end{figure}

The process $y(x)=\int^x_0b_H(s)ds$ in (0,1) was generated with the same
parameters $H, T$ and $N$ as in the bilateral case. The computation
is illustrated, as above, by the distribution of the position
of the maximum $G$ for $y(x)$ with $H=0.1-0.4$ in Fig. 2 (left) and
with $H=0.5-0.9$ in Fig. 2 (right). The plots clearly reveal
the influence of the atom $p_0(T)$
at zero for $H>0.5$, which impedes
to see the exponent  $\theta =\theta (H)$ in the power
law asymptotics of $F_G$  near zero. One can nevertheless assert that
$\theta (H)\to 0$ at the endpoints $H=0$ and 1.
The maximum likelihood estimates of $\theta$ in the interval
of $G:$ $10^{-3}-10^{-2}$ are as follows:

\begin{eqnarray}
\nonumber
H&:&     \quad 0.1 \quad 0.2  \quad 0.3  \quad 0.4  \quad 0.5\\
\nonumber
\hat{\theta}\,\,&:& \quad .09  \quad .13  \quad .18  \quad .21  \quad .23 \\
\nonumber
\theta_0\,&:&  \quad .09  \quad .16  \quad .21  \quad .24  \quad .25
\nonumber
\end{eqnarray}

We also list the hypothetical values of $\theta : \theta_0=H(1-H)$
for comparison  purposes in the above table. The exact result due to
Sinai [18]: $\theta =.25$ for $H=1/2$ shows that we can still take the
error of $\hat{\theta}$ equal to 0.03. In that case the hypothetical
estimates $\tilde{y}(l)$ do not contradict the empirical ones.

Speaking in terms of the process $\tilde{y}(l)$, which can
be obtained by replacing the time in $y(x)$ with local time of
$b_H(x)$, the hypothesis $\theta_0(H)$ means that the position of
the maximum $\tilde{G}_L$ in $\tilde{y}(l)$, $0<l<L$
has the property

\begin{eqnarray}
P(\tilde{G}_L < 1) = L^{-H}{\cal L}_L, \quad L\to \infty
\eqnum{19}
\end{eqnarray}
where $\tilde{{{\cal L}}}_L$ is a slowly increasing function of $L$.
In the limiting case $H\to 0$ the local time $l(x)\simeq x$,
hence $\tilde{y}(l)=\xi l$
where $\xi$ is a Gaussian variate. Consequently,
$P(\tilde{G}_L<1)\sim 1/2$,  $H\to 0$, which is
consistent with the hypothesis (19).

It would be natural to expect an analytical dependence
of $\theta$ on $H$ for the IFBM process. Consequently, the hypothesis
$\theta_0(H)$ can also be extended to cover the case $H>1/2$.
That extrapolation is exact for $H=1$, because $y(x)=\xi x^2/2$,
where $\xi$ is a Gaussian variate, so that $P(G<x)=1/2$
for $H=1$. The last result corresponds to $\theta (1)=0$.

The rigorous result guarantees that $\theta (H)/(1+H)$ is decreasing
for $H>1/2$.

{\bf Statement 4.} (a) The distribution $F_M(x|H)$ for the maximum of

\begin{eqnarray}
y_1(t) = \sqrt{2H+2} \int^t_0b_H(s)ds, \quad 0<t<1
\eqnum{20}
\end{eqnarray}
increases with increasing $H$ in the interval (1/2,\,1)
for any fixed $x>0$.

(b) If $F_M(x^{1+H}|H)=x^{\theta (H)}{\cal L}(x)$, $x\downarrow 0$
where ${\cal L}(x)$ is a slowly varying function, then
$\theta (H)/(1+H)$ decreases with increasing $H$ in (1/2,\,1).

{\it Proof.} The process (20) differs from $y(x)$ by the
normalization $E|y_1(1)|^2=1$. Let $\xi_q(x)=y_1(x^\theta)$ where
$q=2H+2$, $\theta =q_0/q$, $q_0=2H_0+2$ and $H>H_0$.

Since IFBM is self-similar with parameter $h=H+1$, one has
$E|\xi_q(x)|^2 = |x|^{q_0} = E|\xi_{q_0}(x)|^2$

We show in the Appendix that

\begin{eqnarray}
E\xi_q(x)\xi_q(y) \ge E\xi_{q_0}(x)\xi_{q_0}(y)
\eqnum{21}
\end{eqnarray}
when $H>1/2$. In that case the Slepian lemma [13] yields

\begin{eqnarray}
P(\max_{[01]}\xi_q(x)<u) \ge P(\max_{[01]}\xi_{q_0}(x)<u).
\nonumber
\end{eqnarray}
However,
$\max_{[01]}\xi_q(x) = \max_{[01]}\sqrt{2H+2} \int^x_0b_H(s)ds$
which proves the first part of the statement. The second part is
an obvious corollary of the first.

\bigskip

{\bf 5. Conclusion}

\bigskip

We were testing the hypothesis that the maximum $M$ of the integral
of fractional Brownian motion with index $H$ has the distribution
$F_M(x^{(1+H)})=x^{\theta (H)}{\cal L}_H(x)$,  $x\to 0$
in $I\ni 0$, where ${\cal L}_H$ is a slowly varying function.
We have presented theoretical arguments and computational evidence
to support and refine the hypothesis as follows: $\theta (H)=1-H$
for $I=(-1,1)$ and $\theta (H)=H(1-H)$ for $I=(0,1)$.
The computational part of
the problem faces the following difficulties. Due to discrete time
(a step of $\Delta$), the analogue of $M$ in the grid case $(\tilde{M})$
has a nonzero probability $P(\tilde{M}=0)$ of the same order as
$F_M(\Delta^{1+H})$.
It causes a slow (power law) convergence of the
distributions of $\tilde{M}$ and $M$.
Further, there exists an interval
$(0, t_0(\Delta,H))$, $t_0\to 0$ as
$\Delta \to 0$, where $P(\tilde{M}<x^{1+H})$
has a power law behavior with exponent
$\tilde{\theta}(H)>\theta (H)$, $H>1/2$. This makes
the choice of the interval where $\theta (H)$ is to be estimated more
difficult. When $\Delta =0.00012$ and the number of samples is
$N=50 000$, the accuracy of $\theta (H)$ is $\sim 0.03$ for the
most favorable range of $H$: $H\le 1/2$.

\bigskip

{\bf 6. Appendix}

\bigskip

{\it Proof of Theorem 2}. A distribution function (here $F_G$) is
differentiable almost everywhere, because of monotonicity.
Suppose this is true for the point $x_0$. Consider
$x<x_0,\,\lambda =x_0/x >1$.
The self-similarity of $\xi (x)$ yields

\begin{eqnarray}
P(G(0,1)\in dx)=P(G(0,\lambda )\in \lambda dx) \le
P(G(0,1)\in \lambda dx) = \psi (x_0)\frac {x_0}{x} dx
\eqnum{22}
\end{eqnarray}
Here, $G(a,b)$ is the leftmost position of the supremum of
$\{ \xi (x), x\in (a,b)\}$. Consequently, the distribution of $G(0,1)$
is absolutely continuous in $(0,x_0)$.
Points like $x_0$ are dense in (0,1). Consequently,
$F_G(dx)=\psi (x)dx$,\,$x\in (0,1)$.

In virtue of (22), $\psi (x)x$ is a nondecreasing function, i.e.,
the discontinuities in $\psi$  are at most denumerable,
while finite limits on the left and the right exist at the
discontinuity points. The fact that the increments of $\xi (x)$ are
stationary yields

\begin{eqnarray}
\nonumber
P\biggr \{ G(0,1)\in dx\biggl \}&=&P\{ G(-a,\,1-a )\in dx-a\}
\le P\{ G(0,1-a)\in dx-a\} \\
         &=&P\biggr \{ G(0,1)\in d\biggr (\frac {x-a}{1-a}\biggl )\biggl \}
\nonumber
\end{eqnarray}
for any $0<a<1$. One has

\begin{eqnarray}
\psi (x) \le \psi \biggr (\frac {x-a}{1-a}\biggl ) \frac {1}{1-a}
\nonumber
\end{eqnarray}
at continuous points of $\psi$.
Multiply both parts by $(1-x)$:

\begin{eqnarray}
(1-x)\psi (x) \le
\psi \biggr (\frac {x-a}{1-a}\biggl ) \biggr (1-\frac {x-a}{1-a}\biggl )
=\psi (y)(1-y), \quad y=\frac {x-a}{1-a}<x.
\nonumber
\end{eqnarray}
Combining both inequalities, one gets

\begin{eqnarray}
\psi (x) \le \psi (y)\max
\biggr (\frac {y}{x}, \frac {1-y}{1-x}\biggl )
\eqnum{23}
\end{eqnarray}
at all points where $x$ and $y$ are continuous. In particular,
$\psi (x+0)  \le \psi (x-0) \le \psi (x+0)$,
i.e., $\psi$ is continuous in (0,1). If
$\psi (x_0)=0$, $x_0\in (0,1)$, then one has $\psi (x)=0$, $x\in (0,1)$
from (23). Consequently, the following alternative holds: either
$\psi \equiv 0$ or $\psi >0$ in (0,1). The second part of Theorem 2
is an immediate corollary of the first part and the
self-similarity of $\xi (x)$, see [15].

{\it The proof of (21) in Statement 4}.

Let $\xi_q(t)=\sqrt{q}\int^\tau_0b_H(s)ds$,\, $\tau =t^\theta$,
\,$\theta =q_0/q<1$,  $q=2H+2$. The correlation function $\beta_q(t,s)$
of $\xi_q(t)$,  can be written as

\begin{eqnarray}
2\,t^{-q_0}\beta_q(t,s) = \frac{q_0}{q_0-\theta}(\rho^\theta +
\rho^{q_0-\theta}) + [(1-\rho^\theta)^{q_0/\theta} - (1+\rho^{q_0})]
\frac{\theta}{q_0-\theta}
\eqnum{24}
\end{eqnarray}
where $\rho =s/t$. Because $\beta_q$ is symmetric in $t, s$, we put
$\rho \le 1$. We will show that  $\beta_q(t,s)>\beta_{q_0}(t,s)$,
if $q>q_0>3$ or, which amounts to the same
thing, $H>H_0>1/2$.

We have in virtue of (24):

\begin{eqnarray}
2\,t^{-q_0}[\beta_q(t,s) - \beta_{q_0}(t,s)] =
(q-q_0)(1-\rho^\theta)(q_0-1)^{-1}R_1 + q_0(q_0-1)^{-1}R_2.
\nonumber
\end{eqnarray}
Here  \,
$R_1=(1-y^{q-1}-\bar{y}^{q-1})(q-1)^{-1}-(\bar{y}^{q_0-1}-\bar{y}^{q-1})
(q-q_0)^{-1}$,\, $y=\rho^\theta$,\, $\bar{y}=1-y$ and

\begin{eqnarray}
R_2 = \int^1_\theta [\rho^\alpha - \rho^{q_0-\alpha} -
(1-\rho^\alpha)^{q_0-1}\rho^\alpha]\,d\alpha\,\ln 1/\rho.
\nonumber
\end{eqnarray}

We now are going to show that $R_1\ge 0$, and $R_2\ge 0$, if $H\ge 1/2$.

{\it Consider} $R_2$. Put $\rho^\alpha=u$. Since
$0<\theta <\alpha <1$ and $0<\rho <1$, it
follows that $0<u<1$. The integrand in $R_2$ becomes

\begin{eqnarray}
\nonumber
u[1-u^{q_0/\alpha -2}-(1-u)^{q_0-1}] &\ge & u[1-u^{q_0-2}-(1-u)^{q_0-2}] \\
&\ge &u[1 - \max(1, 2^{3-q_0})] \ge 0
\nonumber
\end{eqnarray}

The last estimate is true, because $3-q_0=1-2H\le 0$ when $H>1/2$.
Consequently, $R_2 \ge 0$.

{\it Consider} $R_1$. The function $R_1(y)$ is positive around 0 and 1:

\begin{eqnarray}
R_1 = \left \{ \begin{array}{ll}
y^2(q_0-1)/2 + O(y^{q-1}), \quad  y\to 0 \\
\bar{y} +  O(\bar{y}^{q_0-1}),\qquad \, \qquad \quad \,\,\,\bar{y}\to 0.
\end{array}\right.
\eqnum{25}
\end{eqnarray}

Consequently, $R_1\ge 0$, if the function has a single local extremum
in (0,1). Let $z=(1-y)^{-1}\in (1,\infty)$. Then

\begin{eqnarray}
-z^{q_0-2}\frac{d}{dy}R_1 = [(z-1)^{q-2}-1] - [z^{q-q_0}(q_0-1)-q-1]
(q-q_0)^{-1}:=f(z)
\nonumber
\end{eqnarray}

We now show that $f(z)$ has a single root in $(1,\infty )$. The function

\begin{eqnarray}
f(z) = \left \{ \begin{array}{ll}
-(z-1)(q_0-1) + O((z-1)^{2H}), \quad \,\,\, z\to 1 \\
(z-1)^{q-2}(1+o(1)),\quad \qquad \qquad \qquad z\to \infty
\end{array}\right.
\eqnum{26}
\end{eqnarray}
changes sign in $(1,\infty)$. The equation $f'(z)=0$ or

\begin{eqnarray}
(q-2)(z-1)^{q-3} = (q_0-1)z^{q-q_0-1}
\eqnum{27}
\end{eqnarray}
determines the local extremums of $f(z)$ in $(1,\infty)$. Two strictly
monotone functions occur in (27): on the left is a function that
increases from 0 to $\infty$, because $q-3=2H-1>0$, while the nonnegative
function on the right decreases toward zero at infinity, because
$q-q_0-1=2(H-H_0)-1<2\cdot1/2 - 1<0$.
In that case, however, (27) has the single root $1<z^*<\infty$.
In virtue of (26) $z^*$ is the point of minimum, where $f(z^*)<0$.
The function $f(z)$ is strictly increasing  from $f(z^*)<0$ to $\infty$
in $(z^*,\infty)$, so that the equation $f(z)=0$ has a single root,
as was to be proved.

\bigskip

{\bf Acknowledgements.} This research was supported by
the James McDonnell Foundation within the framework
of the 21st Century Collaborative Action Award for Studying Complex
Systems (project "Understanding and Prediction of Critical Transitions
in Complex Systems"), by the National Science Foundation
(grant EAR 9804859) and in part by the Russian Foundation for
Basic Research (grant 99-01-00314).

\newpage

\bigskip

\centerline{\bf REFERENCES}

\begin{description}

\item
 1.\,\, Ammar G.S., and W.B. Cragg. Superfast solution of real positive
Toeplitz systems. {\it SIAM J. Matrix Annal. Appl.}, {\bf 9}:1, 61-76 (1988).

\item
 2.\,\, Bardet J.M., G. Lang, A. Philippe, and M.S. Taqqu. Generators
of long-range dependent processes:
a survey, in: Donkham P., G. Oppenheim, and M. Taqqu (eds.),
{\it Long-range Dependence: Theory and Applications}, vol. 1, 2002,
579-623, Birkhauser Production.

\item
 3.\,\, Bertoin J. The inviscid Burgers equation with Brownian
initial velocity. {\it Commun. Math. Phys.} {\bf 193}, 397-406 (1998).

\item
 4.\,\, Fernique X. Regularite des trajectoires des fonctions aleatoires
gaussiennes. {\it Lecture Notes in Mathematics}, vol. 1480: 2-187 (1975).

\item
 5.\,\, Handa K. A remark on shocks in inviscid turbulence, in: N.
Fitzmaurice et al. (eds.), {\it Nonlinear Waves and Turbulence},
pp. 339-345, 1993. Birkhauser, Boston.

\item
 6.\,\, Hu X., and S.J. Taylor. The multifractal structure of stable
occupation measure. {\it Stochastic Processes Appl.} {\bf 66},
283-299 (1997).

\item
 7.\,\, Hu X., and S.J. Taylor. The multifractal structure of a general
subordinator. {\it Stochastic Processes and their Appl.} {\bf 88},
245-258 (2000).

\item
 8.\,\, Jaffard S. The multifractal nature of Levy processes. {\it Probab.
Theory Relat. Fields} {\bf 114}:2, 207-227 (1999).

\item
 9.\,\, Isozaki Y. Asymptotic estimates for the distribution of additive
functionals of Brownian motion by the Wiener-Hopf factorization method.
{\it J. Math. Kyoto Univ.} {\bf 36}:1, 211-227 (1996).

\item
10.\,\, Isozaki Y., and S. Kotani. Asymptotic estimates for the first
hitting time of fluctuating additive functionals of Brownian motion.
{\it Lecture Notes in Math.}, 1729, 374-387 (2000).

\item
11.\,\, Isozaki Y., and S. Watanabe, An asymptotic formula for the
Kolmogorov diffusion and a refinement of Sinai's estimates for the
integral of Brownian motion. {\it Proc. Japan Acad.} {\bf 70A}, 271-276 (1994).

\item
12.\,\, Kahane J.-P. {\it Some Random Series of Functions.} 2nd ed.,
Cambridge University Press (1985).

\item
13.\,\,\, Leadbetter M., G. Lindgren, H. Kootzen.
{\it Extremes and related properties of random sequences and processes}.
Springer-Verlag Inc. (Springer ser. in Statistics) (1986).

\item
14.\,\, Majumdar S.N. Persistence in nonequilibrium systems.
{\it Current Science} {\bf 77}:3, 370-375 (1999).

\item
15.\,\, Molchan G. Maximum of a fractional Brownian motion:
probabilities of small values. {\it Commun. Math. Phys.} {\bf 205},
97-111 (1999).

\item
16.\,\, Molchan G.M., and Yu.I. Golosov, Gaussian stationary processes
with asymptotic power spectrum. {\it Soviet Math. Dokl.}
{\bf 10}:1, 134-137 (1969).

\item
17.\,\, She Z., E. Aurell, and U. Frisch. The inviscid Burgers equation with
initial data of Brownian type, {\it Commun. Math. Phys.} {\bf 148},
623-642 (1992).

\item
18.\,\, Sinai Ya.G. Statistics of shocks in solutions of the inviscid Burgers
equation, {\it Commun. Math. Phys.} {\bf 148}, 601-621 (1992).

\item
19.\,\, Vergassola M., B. Dubrulle, U. Frisch, and A. Noullez. Burgers'
equation, Devil's staircases and the mass distribution for large-scale
structures. {\it Astron. Astrophys.} {\bf 289}, 325-356 (1994).

\item
20.\,\, Winkel M. Limit clusters in the inviscid Burgers turbulence with
certain random initial velocities. {\it J. Statistical Phys.}
{\bf 107}, no. 3/4, 893-917 (2002).

\item
21.\,\, Woyczynski W.A. {\it Burgers-KPZ turbulence. G\"ottingen
Lectures}. Springer, (Lectures notes in mathematics; 1700) (1998).

\item
22.\,\, Xiao, Y. H\"{o}lder conditions for the local times and the
Hausdorff measure of the level sets of Gaussian random fields.
{\it Probab. Theory Relat. Fields}, {\bf 109}, 129-157 (1997).

\end{description}

\end{document}